\documentclass[12pt]{article}
\pdfoutput=1
\usepackage{color}
\usepackage{times}
\usepackage{amsmath}
\usepackage{amssymb}
\usepackage{amscd}
\usepackage{url}
\usepackage{doi}

\usepackage{bm}

\title{A Local Inverse Formula and a Factorization}

\author{Gilbert Strang \quad and \quad  Shev MacNamara}

\date{}

\begin{document}




\maketitle

\vspace{1cm}
\begin{center}
With congratulations to Ian Sloan!
\end{center}
\vspace{1cm}

\begin{abstract}
 When a matrix has a banded inverse there is a remarkable formula that quickly computes that inverse, using only local information in the original matrix.
 This local inverse formula holds more generally, for matrices with sparsity patterns that are examples of chordal graphs or perfect eliminators.
 The formula has a long history going back at least as far as the  completion problem for covariance matrices with missing data.
Maximum entropy estimates, log-determinants, rank conditions, the Nullity Theorem and wavelets are all closely related, and the formula has found wide applications in machine learning and graphical models.
 We describe that local inverse and explain how it can be understood as a matrix factorization.
\end{abstract}

\clearpage
\newpage
\section{Introduction}
\label{sec:1}
\vspace{0.5cm}
Here is the key point in two sentences.
If a square matrix $M$ has a tridiagonal inverse, then $M^{-1}$ can be determined from the tridiagonal part $M_0$ of the original $M$.
The formula for $M^{-1}$ is ``local'' and fast to compute --- it uses only $1 \times 1$ and $2 \times 2$ submatrices (assumed invertible) along the main diagonal of $M_0$.

Outside of $M_0$, the entries of $M$ could be initially unknown (``missing data'').
They are determined by the requirement that $M^{-1}$ is tridiagonal.
That requirement maximizes the determinant of the completed matrix $M$: the entropy.

This theory (developed by others) extends to all chordal matrices: the non-zero positions $(i,j)$ in $M_0$ correspond to edges of a chordal graph.
In applications these come primarily from one-dimensional differential and integral equations.
We believe that this special possibility of a local inverse should be more widely appreciated.
It suggests a fast preconditioner for more general problems.


In our first examples of this known (and surprising) formula, $M^{-1}$ will be block tridiagonal:
\vspace{0.5cm}
\begin{equation}
M^{-1} =
\left(
\begin{tabular}{ccc}
$B_{11}$ & $B_{12}$ & $\bm{ 0 }$\\
 $B_{21}$ & $B_{22}$ & $B_{23}$  \\
 $\bm{0}$ & $B_{32}$ & $B_{33}$
\end{tabular}
\right).
\label{eq:block:3x3:matrix:inverse}
\end{equation}
\vspace{0.5cm}

This imposes a strong condition on $M$ itself, which we identify now.
$M$ will be written in the same block form with $n$ square blocks along the main diagonal: $n=3$ above.
When all blocks are $1 \times 1$, the entries of $M^{-1}$ are known to be cofactors of $M$ divided by the determinant of $M$.
Those zero cofactors (away from the three central diagonals) mean that $M$ is ``\textit{semiseparable}:'' all submatrices that don't cross the main diagonal have rank $1$.

In other words, all the $2 \times 2$ submatrices of $M$ (that do not cross the main diagonal) will be singular.
There is a well-developed theory for these important matrices \cite{SemiSeparableBook} that allows wider bands for $M_0$ and $M^{-1}$.

Here are equivalent conditions on $M$ that make $M^{-1}$ \textit{block tridiagonal}.
The key point is that the entries in $M_0$ determine all other entries in $M$.
Those entries are shown explicitly in condition $2$.

\vspace{1cm}
\begin{enumerate}
\item The completion from $M_0$ to $M$ maximizes the determinant of $M$ (the \textit{entropy}).
\item The completion for $n=3$ is given by
\begin{equation}
M=
\left(
\begin{tabular}{ccc}
$M_{11}$ & $M_{12}$ & \; $\bm{ M_{12}M_{22}^{-1}M_{23} }$\\
 $M_{21}$ & $M_{22}$ & $M_{23}$  \\
 $\bm{M_{32} M_{22}^{-1}M_{21}}$ \; & $M_{32}$ & $M_{33}$
\end{tabular}
\right).
\label{eq:block:3x3:matrix:completion}
\end{equation}
Applying this rule recursively outward from the main diagonal, $M$ is obtained from $M_0$ for any matrix size $n$.
\item The completed entries $M_{13}$ and $M_{31}$ minimize the ranks of
\[
\left(
\begin{tabular}{ll}
$M_{12}$ &  $M_{13}$\\
  $M_{22}$ & $M_{23}$
\end{tabular}
\right)
\quad
\textrm{ and }
\quad
\left(
\begin{tabular}{cc}
$M_{21}$ &  $M_{22}$\\
  $M_{31}$ & $M_{32}$
\end{tabular}
\right).
\]
\end{enumerate}
This extends the ``zero cofactor'' condition to the block case.
For ordinary tridiagonal matrices $M$, all $2 \times 2$ blocks  (except those crossing the main diagonal) have rank 1.
The Nullity Theorem \cite{StrangNguyenInterplay2004} says that the dimensions of the null spaces of these $2 \times 2$ block matrices match the number of columns in the zero blocks in $M^{-1}$.



Early motivation for these matrix problems came in statistics where a covariance matrix might have missing entries because  complete data was not available.
For instance in finance, perhaps two assets are not traded frequently enough or not on sufficiently comparable time-scales to provide data that would lead to a sensible estimate of a covariance.
An example of an incomplete covariance matrix could be (although more complicated patterns of missing entries are possible)
\[
M_0= \left(
\begin{tabular}{ccc}
$M_{11}$ & $M_{12}$ & \; $\bm{ ?}$\\
 $M_{21}$ & $M_{22}$ & $M_{23}$  \\
 $\bm{?}$ \; & $M_{32}$ & $M_{33}$
\end{tabular}
\right).
\]
A na\"ive first remedy is to replace the missing entries in $M_0$ by zeros, but that is usually not a good idea; amongst other issues that choice is not guaranteed to always result in a positive definite completion.
Dempster \cite{Dempster72} suggested completing to a covariance matrix by instead inserting zeros in the inverse matrix in positions that correspond to missing values in the original incomplete covariance matrix.
 For this example, that leads to a sparsity pattern for $M^{-1}$ displayed in \eqref{eq:block:3x3:matrix:inverse}, and eventually to the completion to $M$ displayed in \eqref{eq:block:3x3:matrix:completion}.
In general, the entries of the inverse covariance matrix (the \textit{concentration matrix} or the \textit{precision matrix}) can be interpreted as the \textit{information}, so setting these entries to zero reflects the situation that in the absence of data we have no information.
More than that, in the multivariate Gaussian case where a vector $x \in \mathbb{R}^d$ has probability density
\[
p(x) = \frac{1}{\sqrt{(2 \pi)^d} \sqrt{\det M}} \exp(-x^\top M^{-1} x /2),
\]
the entropy $ \int p(x) \log p(x) \mbox{d} x$ (an integral in $d$-dimensions), of the distribution is maximized by maximizing the determinant.
The zeros in the inverse matrix are a consequence of maximizing the determinant \cite{JohnsonProceedings1989} subject to the constraint of being consistent with the initial data.
That seems intuitively satisfying because maximizing an entropy corresponds in some sense to assuming as little as possible while  remaining consistent with the partial data.
This also leads to the maximum likelihood estimator.

Zeros in the concentration matrix $M^{-1}$ correspond to \textit{conditional independence} and the non-zero pattern of the concentration matrix corresponds to edges in the graph of the associated Gaussian Markov random field \cite{SpeedKiiveri1986}.
Extending these ideas to estimate covariance matrices in \textit{high dimensions} is an important and active line of work, connecting to methods that impose sparsity on the concentration matrix via $l_1$-regularization while still optimizing log-determinant objective functions \cite{ravikumar2011high,friedman2008sparse}.

Other references include: Gohberg \textit{et. al.} \cite{DymGohberg1981,EidelmanGohbergHaimovici2013}, Johnson and Lundquist \cite{JohnsonProceedings1989,JohnsonLocalInverse1998}, Lauritzen \cite[page 145]{GraphicalModelsLauritzenBook},  Speed and Kiiveri \cite{SpeedKiiveri1986} and Strang and Nguyen \cite{StrangNguyenInterplay2004,Strang:2010aa}.


\section{The Local Inverse Formula}\label{sec:2}
At first sight it is hard to believe that the inverse of an $n \times n$ matrix (or block matrix) can be found from ``local inverses.''
But if $M^{-1}$ is a tridiagonal (or block tridiagonal) matrix, that statement is true.
The only inverses you need are $1 \times 1$ and $2 \times 2$ along the main diagonal of $M$.
The $1 \times 1$ inverses, $M_{2,2}^{-1}$, $\ldots$, $M_{n-1,n-1}^{-1}$, come from the diagonal blocks.
 The $2 \times 2$ inverses $Z_i^{-1}$ come from adjacent blocks:
\[
Z_{i}^{-1} =
\left(
\begin{array}{cc}
M_{i,i} & M_{i,i+1}\\
M_{i+1,i} & M_{i+1,i+1}
\end{array}
\right)^{-1}.
\]
In other words, we only need the tridiagonal part of $M$  to find the tridiagonal matrix $M^{-1}$.

From $n-2$ inverses $M_{i,i}^{-1}$ and $n-1$ inverses $Z_{i}^{-1}$, here is the \textbf{local inverse formula} (when $n=3$):
\begin{eqnarray}
M^{-1} &=&
\left(
\! \! \! \! \!
\begin{array}{ccc}
\left(
\! \! \! \!
\begin{tabular}{ccc}
$M_{11}$ & $M_{12}$ & \\
 $M_{21}$ & $M_{22}$ &
\end{tabular}
\! \! \! \! \! \! \! \! \! \!
\right)^{-1}
& \\
& \\
\end{array}
\right)
+
\left(
\begin{array}{crr}
&\\
&  \left(
\! \! \! \! \! \! \! \! \!
\begin{tabular}{crr}
  & $M_{22}$ & $M_{23}$  \\
  & $M_{32}$ & $M_{33}$
\end{tabular}
\! \! \! \!
\right)^{-1}
\end{array}
\! \! \! \! \!
\right)
-
\left(
\! \! \! \! \!
\begin{tabular}{ccc}
\\
& $M_{22}^{-1}$ &  \\
&   &
\end{tabular}
\! \! \! \! \!
\right)  \label{eq:block:3x3:matrix:local:inverse:formula}  \\
&=&
\left(
\! \! \! \! \!
\begin{array}{ccc}
\left(
\begin{tabular}{ccc}
 &  & \\
  & $Z_1^{-1}$ &  \\
  &
\end{tabular}
\right)
& \\
& \\
&
\end{array}
\right)
+
\left(
\begin{array}{ccr}
&\\
& & \\
&
\left(
\begin{tabular}{ccc}
  &  &   \\
  &  & $Z_2^{-1}$ \\
  &
\end{tabular}
\right)
\end{array}
\! \! \! \! \! \! \!
\right)
-
\left(
\begin{tabular}{ccc}
& \\
& \\
& $M_{22}^{-1}$ &  \\
&   & \\
&
\end{tabular}
\right).
\nonumber
\end{eqnarray}
We emphasize that $M$ itself need not be tridiagonal.
It rarely is.
The construction of $M$ does \textit{start} with a tridiagonal matrix, $M_0$.
That matrix is completed to $M$ in such a way that $M^{-1}$ is tridiagonal.
It becomes reasonable to expect that $M^{-1}$ depends only on the starting tridiagonal matrix $M_0$.
But still the simplicity of the local inverse formula is unexpected and attractive.

This local formula for $M^{-1}$ can be established in several ways.
Direct matrix multiplication will certainly succeed.
Johnson and Lundquist  \cite{JohnsonLocalInverse1998} show how this $3 \times 3$ block case extends by iteration to larger matrices (with wider bands or general chordal structures, described next).
The present paper looks at the triangular $LDU$ factorization --- which produces banded or chordal factors.
And we view the matrix algebra in the $A^\top C A$ framework that is fundamental to applied mathematics.

The generalisation to matrices $M^{-1}$ with five non-zero block diagonals is straightforward.
Thus $M_0$ is pentadiagonal, and its extension to $M$ is determined so that $M^{-1}$ is also pentadiagonal.
Then the local inverse formula goes directly from $M_0$ to $M^{-1}$, bypassing the completed matrix $M$.
The formula involves the $2 \times 2$ inverses $Z_i^{-1}$ together with the $3 \times 3 $ inverses $Y_i^{-1}$.
The submatrices $Y_i$ come from three adjacent rows and columns ($i,i+1,i+2$) of $M_0$ and $M$.
The local inverse formula assembles the inverse as before (displayed for $n=4$):
\begin{eqnarray}
M^{-1} &=&
\left(
\! \! \! \! \!
\begin{array}{ccc}
\left(
\begin{tabular}{ccc}
 &  & \\
  & $Y_1^{-1}$ &  \\
  &
\end{tabular}
\right)
& \\
& \\
& \\
\end{array}
\right)
+
\left(
\begin{array}{ccr}
&\\
& \\
\left(
\begin{tabular}{ccc}
  &  &   \\
  &  & $Y_2^{-1}$ \\
  &
\end{tabular}
\right)
\end{array}
\! \! \! \! \! \! \!
\right)
-
\left(
\begin{tabular}{ccc}
& \\
& \\
& $\left( Z_{2}^{-1} \right)$ &  \\
&  \\
&
\end{tabular}
\right) \nonumber \\
&=& \textrm{5-diagonal matrix.}
\nonumber
\end{eqnarray}
The formula extends to wider bands in $M_0$ in a natural way.
Beyond that come `staircase matrices' that are unions of overlapping square submatrices $Y_i$ centered on the main diagonal.
The sizes of the $Y_i$ can vary and the overlaps (intersections) are the $Z_i$.
The  inverse formula remains correct.

The ultimate extension is to \textit{chordal matrices} $M_0$ and $M^{-1}$ \cite{JohnsonLocalInverse1998}.
Their non-zero entries produce a \textit{chordal graph} \cite{BlairPeytonReport1993,PeterBartlett2009,KollerFriedman2009}.
Beyond that we cannot go.
Two equivalent definitions of the class of chordal matrices are:
\begin{itemize}
\item Suppose $M_0$ has non-zero entries in positions $(i_0,i_1), (i_1,i_2), \ldots, (i_m,i_0)$.
If $m \ge 4$ then that closed path has a ``shortcut'' chord from an $i_J$ to an $i_L \ne i_{J+1}$ for which $M_0(i_J,i_L) \ne 0$.
\item There are permutations $P$ and $Q^\top$ of the rows and columns of $M_0$ so that the matrix $A= P M_0 Q^{\top}$ allows ``\textit{perfect elimination with no fill-in}:''
\[
A = LDU = \textrm{(lower triangular) (diagonal) (upper triangular) }
\]
\[
\textrm{with } L_{ij} = 0 \textrm{ and } U_{ij}=0 \textrm{ whenever } A_{ij}=0.
\]
\end{itemize}
We may assume  \cite{DonRose70}, that $M_0$ comes in this perfect elimination order.
Then it is completed to $M$ in such a way that $M$ has the same elimination order as $M_0$.

\section{Completion of $M$ and triangular factorizations}\label{sec:3}
When does a $3 \times 3$ block matrix $M$ have a tridiagonal inverse?
If the tridiagonal part of $M$ itself is prescribed, the entries in the upper right and lower left corners are determined by the requirement that the corresponding entries in $M^{-1}$ are zero:
\begin{equation}
M =
\left(
\begin{tabular}{ccc}
$M_{11}$ & $M_{12}$ & \; $\bm{ M_{12}M_{22}^{-1}M_{23} }$\\
 $M_{21}$ & $M_{22}$ & $M_{23}$  \\
 $\bm{M_{32} M_{22}^{-1}M_{21}}$ \; & $M_{32}$ & $M_{33}$
\end{tabular}
\right).
\label{eq:block:3x3:matrix}
\end{equation}
It is this completed matrix $M$ (also in \eqref{eq:block:3x3:matrix:completion}) that multiplies the matrix in \eqref{eq:block:3x3:matrix:local:inverse:formula} to give the identity matrix and verify the local inverse formula.
Suppose $M$ is block upper triangular: call it $U$, with unit diagonal blocks.
Then the matrices and the local inverse formula become particularly simple.
Here are the incomplete $U_0$, the completed $U$ and the inverse $U^{-1}$:
\[
U_{0} =
\left(
\begin{tabular}{ccc}
$I$ & $U_{12}$ & $\bm{?}$\\
 $0$ & $I$ & $U_{23}$  \\
 $0$ & $0$ & $I$
\end{tabular}
\right)
\]
\begin{equation}
U =
\left(
\begin{tabular}{ccc}
$I$ & $U_{12}$ & $\bm{U_{12} U_{23}} $\\
 $0$ & $I$ & $U_{23}$  \\
 $0$ & $0$ & $I$
\end{tabular}
\right)
\label{eq:U}
\end{equation}
\[
U^{-1} =
\left(
\begin{tabular}{ccc}
$I$ & $-U_{12}$ & $\bm{0}$\\
 $0$ & $I$ & $-U_{23}$  \\
 $0$ & $0$ & $I$
\end{tabular}
\right).
\]
The local inverse formula separates $U^{-1}$ in three parts:
\begin{equation}
U^{-1}=
\left(
\begin{tabular}{ccc}
$I$ & $-U_{12}$ & $0$\\
 $0$ & $I$ & $0$  \\
 $0$ & $0$ & $0$
\end{tabular}
\right)
+
\left(
\begin{array}{cccc}
0 & &0 & 0\\
 0 & &I & -U_{23}  \\
 0 & &0 & I
\end{array}
\! \! \! \!
\right)
-
\left(
\begin{tabular}{ccccc}
$0$ & & $0$ & & $0$\\
 $0$ & & $I$ & & $0$  \\
 $0$ & & $0$ & & $0$
\end{tabular}
\right).
\label{eq:local:inverse:U}
\end{equation}
We vainly hoped that this simple idea could apply to each factor of  $M=LDU$ and produce factors of $M^{-1}$.
That idea was destined to fail --- the correct factors mix upper with lower (just as elimination does).
Still it would be attractive to understand the general chordal case through its triangular factors.
The key property of ``no fill-in'' distinguishes chordal matrices in such a beautiful way.

\section*{Example}
Consider the $3 \times 3$ matrix
\begin{equation}
M_0 \equiv
\frac{1}{4}
\left(
\begin{tabular}{ccc}
3 & 2 & \textbf{?}    \\
2 & 4 & 2     \\
\textbf{?} & 2 & 3
\end{tabular}
\right)
\nonumber
\end{equation}
that is completed to
\begin{equation}
M =
\frac{1}{4}
\left(
\begin{tabular}{ccc}
3 & 2 & 1    \\
2 & 4 & 2     \\
1 & 2 & 3
\end{tabular}
\right),
\nonumber
\label{eq:simple:3x3:example}
\end{equation}
with inverse
\begin{equation}
M^{-1} =
\left(
\begin{tabular}{ccc}
2 & -1 & 0    \\
-1 & 2 & -1     \\
0 & -1 &2
\end{tabular}
\right).
\nonumber
\end{equation}
In this symmetric example, $U=L^\top$ and $M = L D U = L D L^\top$ where
\begin{equation}
L =
\left(
\begin{tabular}{rrr}
    1 &         0     &    0 \\
  $\frac{2}{3}$  &  1 &         0 \\
           $\frac{1}{3}$  &  $\frac{1}{2}$   & 1
\end{tabular}
\right)
\qquad \textrm{and} \qquad
D = \frac{1}{12}
\left(
\begin{tabular}{rrr}
    9 &         0     &    0 \\
   0  &  8 &         0 \\
         0  & 0   & 6
\end{tabular}
\right) .
\nonumber
\end{equation}
Notice these examples of $L$ and of $U$ satisfy the formats displayed in \eqref{eq:U} and in \eqref{eq:local:inverse:U} (our example illustrates these formats in the scalar case but those formats remain true in the block matrix case).

\vspace{2cm}
\section{The $A^{\top}CA$ framework: A matrix factorization restatement}

Applied mathematics is a broad subject far too diverse to be summarized by merely one equation.
Nevertheless,  $A^\top \, C \,  A$ offers a matrix framework to understand a great many of the classical topics, including: least squares and projections, positive definite matrices and the Singular Value Decomposition, Laplace's equation with $A=\textrm{grad}$ and the Laplacian as $ \nabla^2 = \textrm{div} (\textrm{grad}) = A^\top A $, networks and graph Laplacians \cite{Str16}.
It is therefore satisfying to place the local inverse formula in this framework.

Consider a square $n \times n$ invertible matrix $M$ such that the inverse matrix satisfies the `local inverse formula.'
We will express the local inverse formula as the following factorization of the inverse matrix
\begin{equation}
\label{eq:inverse:factorisation}
M^{-1}   = A^\top \, C^{-1} \, A
\end{equation}
and factorization of the original matrix
\begin{equation}
\label{eq:forward:factorisation}
M  = G^\top \, C \, G .
\end{equation}

Such factorizations are often represented by commutative diagrams.
Here we represent $x = M^{-1} b$  as
\[
\begin{CD}
\mathbb{R}^{m}     @<C^{-1}<<   \mathbb{R}^{m} \\
@V  A^\top  VV        @AAA A \\
\mathbb{R}^{n}     @<M^{-1}<<  \mathbb{R}^{n}
\end{CD}
\]
and we reverse the directions of all four arrows to represent $M x = b$  as
\[
\begin{CD}
\mathbb{R}^{m}     @>C>>   \mathbb{R}^{m} \\
@A  G  AA        @VV G^\top V\\
\mathbb{R}^{n}      @>M>>  \mathbb{R}^{n}
\end{CD}
\]
Notice that in this approach we \textit{start with the inverse matrix} $M^{-1}$, and then we \textit{invert the inverse} to arrive at the original matrix: $(M^{-1})^{-1}=M$.
To describe the factorizations we must identify the matrices $C$, $A$ and $G$, but first we introduce notation.

Our setting is that the non-zero sparsity pattern of $M^{-1}$  is a chordal graph on $n$ nodes, with a clique tree (sometimes called a \textit{junction tree}) on $c_b$ nodes that represent the $c_b$ \textit{maximal cliques} (square submatrices of $M_0$ with no missing entries).
There are $c_b$ `blocks' and $c_o$ `overlaps' in the corresponding local inverse formula.
 Let $c=c_b+c_o$ be the sum of these counts.
Denote these block matrices by $C_k$ for $k=1, \ldots, c$.
Order these matrices so that all the $c_b$ blocks that correspond to maximal cliques come first, and all the $c-c_b$ blocks that correspond to overlaps come last.
Let $d_k$ denote the size of clique $k$ so that $C_k$ is a $d_k \times d_k$ matrix, and let
\begin{equation}
m = \sum_{k=1}^{c} d_k = \sum_{k=1}^{c_b} d_k  \; + \sum_{k=c_b +1}^{c} d_k.
\nonumber
\end{equation}
Note that
$
m>n.
$

Define the $m \times m$ block diagonal matrix
\begin{equation}
C \equiv
\left(
\begin{tabular}{cccclrr}
$C_1$ & &  & &   \\
 & $C_2$ &  &  &    \\
 &  & $\ddots$ & \\
 &  & &   $C_{c_b}$ \\
 \\
 &  & &  & $-C_{c_b+1}$ \\
 \\
  &  & &  &   & $\ddots$ \\
    &  & &  &     & & $-C_{c}$
\end{tabular}
\right).
\nonumber
\end{equation}
The minus signs in front of the blocks in the bottom right corner correspond to overlaps.

Because each block $C_k$ corresponds to a subset of nodes in the original graph, each row ($i=1,\ldots,m$) of $C$ corresponds to a node ($j=1,\ldots,n$) in the original graph.
Define the  $m \times n$ matrix $A$ of $0s$ and $1s$ to encode this correspondence:
\begin{equation}
A_{i,j} \equiv \begin{cases} 1 &\mbox{if node j corresponds to row i of $C$}  \\
0  & \mbox{otherwise}. \end{cases}
\nonumber
\end{equation}
Note that each row of $A$ contains precisely one non-zero entry and that entry is $1$.
The total number of non-zero entries in $A$ is $m$.
Each column of $A$ contains one or more $1s$.

It is a necessary condition for the local inverse formula to apply that all of the blocks $C_k$ be separately invertible.
 Then $C$ is an invertible matrix,  and $C^{-1}$ is the block diagonal matrix with blocks $C_1^{-1}, \ldots, C_c^{-1}$.
With these definitions, the factorization $M^{-1}   = A^\top \, C^{-1} \, A$ in \eqref{eq:inverse:factorisation} is simply matrix notation for the local inverse formula: $M^{-1}$ is ``the sum of the inverses of the blocks, minus the inverses of the overlaps.''

It remains to describe the factorization $M = G^\top C G$ in \eqref{eq:forward:factorisation}.
Intuitively, this is arrived at by reversing the directions of the arrows in the commutative diagram for $M^{-1}$.
It is easy to see that replacing $C^{-1}$ by $C$ will reverse the direction of the arrow at the top of the diagram in a way that correctly inverts the action of $C^{-1}$.

It is not so easy to see that we can find matrices $G^\top$ and $G$ such that the directions of the arrows corresponding to $A$ and to $A^\top$ are reversed with the desired effect.
 Indeed, at first glance that seems to be tantamount to finding the `inverse' of the $A$ matrix, but that is impossible because $A:\mathbb{R}^{n} \rightarrow \mathbb{R}^{m}$ is not a square matrix.
However, there is  redundancy in the action of $A$.
Although $A$ maps from a smaller $n-$dimensional space to a larger $m-$dimensional space, the matrix only has $n$ columns, so the column space reached by $A$ is only an $n-$dimensional subspace of $\mathbb{R}^{m}$.
(Columns of $A$ are independent because each row contains precisely one $1$.)
This makes it possible to choose $G^\top$ so that we only `invert' on the subspace that we need to.
A possible choice is the pseudoinverse of $A$, i.e.
\begin{equation}
\label{eq:F}
F  \equiv (A^\top A)^{-1} A^\top.
\nonumber
\end{equation}
Note that $FA = I_n$ is the $n \times n$ identity matrix (but $FA \ne I_m$).
So $F$ is a left inverse for $A$.
Instead of $F$, we could choose another left inverse of $A$, namely $G^\top$, where $G$ is the matrix described next.

The last step is to now find the matrix $G$ that will `undo' the effect of $A^\top$.
Note that  in our factorization, the matrix $A^\top$ acts only on the range of $C^{-1}A$ (and not on all of $\mathbb{R}^{m}$).
In other words, in our factorization, it is the column space of $C^{-1} A$ that is the `input space' to $A^\top$.
So we only need to invert on that subspace, by
\begin{equation}
\label{eq:G}
G  \equiv C^{-1} A M.
\nonumber
\end{equation}
This choice makes it clear that $G$ has two desirable properties:
\begin{itemize}
\item  the columns of $G$ are linear combinations of the columns of $C^{-1} A$, so the range of $G$ is in the $n-$dimensional subspace of $\mathbb{R}^{m}$ that is reached by $C^{-1} A$, and
\item  $G^{\top} C G = G^\top C (C^{-1} A)M = (G^{\top}A) M = (I) M= M$.
\end{itemize}
The second property, that $L_A C G=M$, is not unique to our choice of $G$ --- it holds for any matrix $L_A$ that is a left inverse of $A$.
(Then we have $L_A C G = L_A C (C^{-1} A)M = (L_AA) M = (I) M= M$.)
We have already seen that $F$ is a left inverse of $A$ so $F$ is a possible choice for a factorization to  recover the original matrix, i.e. $FCG=M$.
To see that $G^\top$ is also a left inverse of $A$,
recall the definition  $G \equiv C^{-1} A M$.
By the rule for a transpose of a product,  $G^\top = M^\top A^\top (C^{-1})^{\top} $.
So
\[
G^\top A =M^\top A^\top (C^{-1})^{\top} A  = M^\top (A^\top C^{-1} A)^\top  = M^\top (M^{-1})^\top = I,
\]
as required.

These choices also have: $A^{\top}G = (A^\top C^{-1} A) M  = M^{-1}  M = I$, and $FG= (A^\top A)^{-1}$.
 And  $AF$ is a projection matrix.

\section*{Example}
We now exhibit the $A^\top C A$ factorization for the same $3 \times 3$ matrix example that we used earlier in \eqref{eq:simple:3x3:example} to demonstrate the $LDU$ factorization
\begin{equation*}
M \equiv
\frac{1}{4}
\left(
\begin{tabular}{ccc}
3 & 2 & 1    \\
2 & 4 & 2     \\
1 & 2 & 3
\end{tabular}
\right) \quad \textrm{with inverse} \qquad
M^{-1} =
\left(
\begin{tabular}{ccc}
2 & -1 & 0    \\
-1 & 2 & -1     \\
0 & -1 &2
\end{tabular}
\right).
\end{equation*}
The graph of non-zeros is a line of three nodes.
The maximal cliques are
\begin{equation}
C_1 \equiv
\frac{1}{4}
\left(
\begin{tabular}{cc}
3 & 2 \\
2 & 4
\end{tabular}
\right) \quad \textrm{and} \quad
C_2 \equiv
\frac{1}{4}
\left(
\begin{tabular}{cc}
 4 & 2     \\
 2 & 3
\end{tabular}
\right).
\nonumber
\end{equation}
The overlap is
\begin{equation}
C_3 \equiv
\frac{1}{4}
\left(
\begin{tabular}{c}
4
\end{tabular}
\right) = (1),
\nonumber
\end{equation}
so in this example $m=2+2+1 = 5$.
The $m \times m $ block diagonal matrix $C$ is
\begin{equation}
C \equiv
\frac{1}{4}
\left(
\begin{array}{ccccr}
3 & 2 & 0  &0 & 0  \\
2 & 4 & 0  &0 & 0  \\
0 & 0 &  4 &2 & 0  \\
0 & 0 & 2  & 3& 0  \\
0 & 0 & 0  &0 & -4
\end{array}
\right).
\nonumber
\end{equation}
The matrix that sends `node space' (the three columns could correspond to the three nodes) to `clique space' (rows $1,2,3,4,5$ correspond to nodes $1,2,2,3,2$) is
\begin{equation}
A \equiv
\left(
\begin{tabular}{ccc}
1 & 0 & 0    \\
0 & 1 & 0    \\
0 & 1 & 0    \\
0 & 0 & 1    \\
0 & 1 & 0    \\
\end{tabular}
\right).
\nonumber
\end{equation}
Direct matrix multiplication confirms that $A^\top C^{-1} A$ does indeed give $M^{-1}$, as expected from the local inverse formula.
In this example
\begin{equation}
F \equiv (A^\top A)^{-1} A^\top =
\left(
\begin{tabular}{ccccc}
1 & 0 & 0 & 0   &0 \\
0 & $\frac{1}{3}$ & $\frac{1}{3}$ & 0   &$\frac{1}{3}$\\
0 & 0 & 0 & 1   &0 \\
\end{tabular}
\right)
\nonumber
\end{equation}
and
\begin{equation}
G \equiv C^{-1}AM =
\left(
\begin{tabular}{rrr}
1 & 0 & 0    \\
0 & 1 & $\frac{1}{2}$    \\
$\frac{1}{2}$  & 1 & 0    \\
0 & 0 & 1    \\
$-\frac{1}{2}$  & $-1$ & $-\frac{1}{2}$
\end{tabular}
\right)
\nonumber
\end{equation}
and  $FCG=G^\top CG=M$.
Typically and in this example, the local inverse formula only applies in going from $M$ to $M^{-1}$, and
\[
M \ne A^\top \, C \, A  =
\frac{1}{4}
\left(
\begin{tabular}{ccc}
3 & 2 & 0    \\
2 & 4 & 2     \\
0 & 2 & 3
\end{tabular}
\right).
\]

\section{Applications}
We now showcase by example some of the especially elegant applications of the local inverse formula.

\section*{Example: A Toeplitz matrix}
Complete the missing entries in $M_0$ to arrive at a first example via
\[
M_0=
\left(
\begin{tabular}{rrrr}
$2$ & $-1$ & {\color{blue} ?} & {\color{blue} ?}\\
 $-1$ & $2$ & $-1$ & {\color{blue} ?} \\
 {\color{blue} ?} & $-1$ & $2$ & $-1$ \\
 {\color{blue} ?} & {\color{blue} ?} &   $-1$ & $2$
\end{tabular}
\right)
\quad
\longrightarrow
\qquad
\left(
\begin{tabular}{rrrr}
$2$ & $-1$ & {\color{blue} $\frac{1}{2}$}  & {\color{blue} $-\frac{1}{4}$}\\
 $-1$ & $2$ & $-1$ & {\color{blue} $\frac{1}{2}$} \\
 {\color{blue} $\frac{1}{2}$} & $-1$ & $2$ & $-1$ \\
 {\color{blue} $-\frac{1}{4}$} & {\color{blue} $\frac{1}{2}$} &   $-1$ & $2$
\end{tabular}
\right) = M
\]
so that the completed matrix has an inverse with zeros in the locations where entries were missing in the original matrix:
\begin{eqnarray}
 M^{-1} &=&
\frac{1}{6}
\left(
\begin{tabular}{rrrr}
$4$ & $2$ & {\color{blue} $0$}  & {\color{blue} $0$}\\
 $2$ & $5$ & $2$ & {\color{blue} $0$} \\
 {\color{blue} $0$} & $2$ & $5$ & $2$ \\
 {\color{blue} $0$} & {\color{blue} $0$} &   $2$ & $4$
\end{tabular}
\right) 
\label{eq:inverse:small:example:with:zeros}  \\
&=&
\frac{1}{6}
\left(
\begin{tabular}{rrrr}
$4$ & $2$ & & \\
 $2$ & $4$  & & \\
& & &  \\
\end{tabular}
\right)
+\frac{1}{6}
\left(
\begin{tabular}{rrrr}
 & & & \\
 & $4$ & $2$ &  \\
  & $2$ & $4$ & \\
  &
\end{tabular}
\right) +
\frac{1}{6}
\left(
\begin{tabular}{rrrr}
 & & & \\
 & & & \\
  & & $4$ & $2$ \\
 &  &   $2$ & $4$
\end{tabular}
\!\! \!
\right) \nonumber \\
& &
\qquad  \qquad \qquad  \qquad \quad
-
\left(
\begin{tabular}{rrrr}
 & & & \\
 & $\frac{1}{2}$ &  & \\
  & & &  \\
 &  &   &
\end{tabular}
\right)
\;\;\;-
\left(
\begin{tabular}{rrrr}
 & & & \\
 &  &  & \\
  & & $\frac{1}{2}$ &  \\
 &  &   &
\end{tabular}
\right).
\label{eq:inverse:small:example:with:zeros:LIF}
\end{eqnarray}
The local inverse formula assembles $M^{-1}$ in \eqref{eq:inverse:small:example:with:zeros:LIF} from the inverses of the three repeating blocks in $M$
\[
\left(
\begin{tabular}{rr}
$2$ & $-1$  \\
 $-1$ & $2$
\end{tabular}
\right)^{-1}
=
\frac{1}{6}
\left(
\begin{tabular}{rr}
$4$ & $2$  \\
 $2$ & $4$
\end{tabular}
\right)
\]
and subtracting the inverses, $(2)^{-1}= 1/2$, of the two overlaps.

To appreciate the significance of those zeros in $M^{-1}$ in \eqref{eq:inverse:small:example:with:zeros}, it helps to recall that the derivative of the determinant with respect to the entries of the matrix is given by a cofactor up to scaling by the determinant (this result comes quickly from the cofactor expansion of the determinant along one row of the matrix, for example).
This leads to an especially simple form of derivative of the \textit{log-determinant}, which in the symmetric case is simply the corresponding entry of the inverse matrix:
\[
\frac{\partial}{\partial a_{ij}} \log \det M = (M^{-1})_{ij}.
\]
Hence zeros in the inverse matrix, such as appear in \eqref{eq:inverse:small:example:with:zeros}, correspond to setting derivatives to zero, which corresponds to a local optima.
The log-determinant is convex on the cone of symmetric positive definite matrices so a local optima is also a global maximum in this case.

This first example suggests a second example, by generalizing to a doubly infinite Toeplitz matrix \cite{DymGohberg1981}.
A \textit{Toeplitz matrix} is constant along diagonals: the $(i,j)$ entry is a function of $(i-j)$, so specifying one row of the matrix completely specifies all entries of the matrix.
In the doubly infinite Toeplitz case, the entries of a row are the Fourier series of an associated  function $s$ known as the \textit{symbol} of the matrix.
 The matrix completion problem becomes a problem of Fourier series for functions.
We must complete the missing Fourier coefficients for a function $s$ so that the Fourier series of the reciprocal function $1/s$ has zero coefficients corresponding to missing entries in the Fourier series  of $s$.
For example,
\[
\! \! \! \! \! \! \! \! \! \! \! \! \! \! \! \! \! \! \! \! \! \! \! \! \! \! \! \!  \! \! \! \! \!
 \! \! \! \! \!  \! \! \! \! \!  \! \! \! \! \!  \! \! \! \! \! \! \! \! \!
 \left(
\begin{tabular}{rrrrrrrrr}
& & &  $\ddots$ \\
$\cdots$ &  {\color{blue} ?} & {\color{blue} ?} & $-1$ & $2$ & $-1$ & {\color{blue} ?}  & {\color{blue} ?} & $\cdots$ \\
& & & & & $\ddots$
\end{tabular}
\right)^{-1} =
\]
\vspace{-0.3cm}
\[
\qquad \qquad \qquad \qquad \qquad \qquad \qquad \,  \left(
\begin{tabular}{ccccccccc}
& & &  $\ddots$ \\
$\cdots$ &  0 & 0& ${\color{blue} ?}$ & ${\color{blue} ?}$ & ${\color{blue} ?}$ & 0& 0 & $\cdots$ \\
& & & & & $\ddots$
\end{tabular}
\right)
\]
\vspace{-0.2cm}
is completed to
\[
\! \! \! \! \! \! \! \! \! \! \! \! \! \! \! \! \! \! \! \! \! \! \! \! \! \! \! \!  \! \! \! \! \!
 \! \! \! \! \!  \! \! \! \! \!  \! \! \! \! \!  \! \! \! \! \! \! \! \! \!
\left(
\begin{tabular}{ccccccccc}
& & &  $\ddots$ \\
$\cdots$ &  {\color{blue} $-\frac{1}{4}$} & {\color{blue} $\frac{1}{2}$} & $-1$ & $2$ & $-1$ & {\color{blue} $\frac{1}{2}$}  & {\color{blue} $-\frac{1}{4}$} & $\cdots$ \\
& & & & & $\ddots$
\end{tabular}
\right)^{-1}
=
\]
\vspace{-0.3cm}
\[
\qquad \qquad \qquad \qquad \qquad \qquad \qquad \,
\frac{1}{6}
\left(
\begin{tabular}{ccccccccc}
& & &  $\ddots$ \\
$\cdots$ &  0 & 0 & {\color{blue} $2$} & {\color{blue} $5$} &{\color{blue} $2$} & 0 & 0& $\cdots$ \\
& & & & & $\ddots$
\end{tabular}
\right).
\]

 The general principle is to complete the symbol
\[
s(x) = \sum_{-\infty}^{\infty} a_k e^{ikx} \qquad \quad \textrm{with inverse} \qquad \quad \frac{1}{s(x)} = \sum_{-\infty}^{\infty} b_k e^{ikx}
\]
so that $\mathbf{b_k=0}$ \textbf{when} $\mathbf{a_k}$ \textbf{was not specified.}
This maximizes the log-determinant
\[
\int_0^{2\pi} \log \sum_{-\infty}^{\infty} a_k e^{ikx} \mbox{d} x
\]
amongst symmetric positive definite Toeplitz matrices.

\vspace{1cm}

\section*{Banded matrices with banded inverse}
In very exceptional cases \cite{Strang:2010aa} a banded matrix can have a banded inverse.
Then the local inverse formula applies in  `both directions' (leading to a class of `chordal matrices with chordal inverse').
This will give a (new?) algorithm for the analysis and synthesis steps in a discrete wavelet transform (known as a filter bank) \cite{DaubechiesBook,MallatWaveletBook,StrangWaveletBook}.
Here is an example of one of the celebrated Daubechies wavelets in this framework.

\section*{Example: a Daubechies wavelet}
Set
\[
s = \sqrt{3},
\quad
B_1 =
\left(
\begin{tabular}{rcr}
$1+s$ & & $3+s$ \\
     $-1+s$ & & $3-s$
\end{tabular}
\right),
\; \textrm{and} \quad
B_2=
\left(
\begin{tabular}{rcr}
$3-s$ & & $1-s$ \\
     $-3-s$ & & $1+s$
\end{tabular}
\right).
\]
Notice $B_1$ and $B_2$ are singular.
Set
\[
t' =  \left(
\begin{tabular}{cccccc}
$-(3+s)$ &$1+s$ & 0 &0 &0 &0
\end{tabular}
\right)
\quad
\textrm{and}
\quad
t =  \sqrt{32} \frac{t'}{||t'||_2}
\]
and
\[
b' =  \left(
\begin{tabular}{cccccc}
0 &0& 0& 0&    $(1+s)$ & $3+s$
\end{tabular}
\right)
\quad
\textrm{and}
\quad
b =  \sqrt{32} \frac{b'}{||b'||_2}.
\]

With these definitions a matrix corresponding to a Daubechies wavelet is
\vspace{0.5cm}
\begin{eqnarray*}
M &=&
\frac{1}{\sqrt{32}}
\left(
\begin{tabular}{ccc}
& $t$ \\
$B_1$ &$B_2$ & $\bm{0}$ \\
$\bm{0}$ & $B_1$& $B_2$\\
& $b$
\end{tabular}
\right) \nonumber \\
&=&
\left(
\begin{tabular}{cccccc}
   -0.8660 &   0.5000 &        0  &       0   &      0    &     0 \\
    0.4830  &  0.8365 &   0.2241 &  -0.1294  &       0   &      0\\
    0.1294 &   0.2241 &  -0.8365  &  0.4830  &       0  &       0\\
         0     &    0  &  0.4830 &   0.8365  &  0.2241 &  -0.1294\\
         0     &    0  &  0.1294  &  0.2241 &  -0.8365 &   0.4830\\
         0     &    0   &      0      &   0  &  0.5000  &  0.8660
\end{tabular}
\right).
\end{eqnarray*}
\vspace{0.5cm}

\noindent  Then, as desired for a wavelet basis, $M$ is orthogonal so $M^{-1}=M^\top$ is also banded.
(There are also important non-orthogonal wavelets with banded $M$ and $M^{-1}$.)

Early motivation for the local inverse formula came from problems with covariance matrices, which are symmetric positive definite.
But the local inverse formula can also apply to matrices that are not symmetric positive definite, as in this Daubechies wavelet matrix example.


More interestingly in the context of our present article, in this example, \textit{the local inverse formula applies in both directions}.
We have
\begin{equation}
A^{\top} C^{-1} A = M^{-1}
\label{eq:local:inverse:formula:forward}
\end{equation}
(this is the local inverse formula that we have come to expect when $M^{-1}$ is chordal) \textit{and}
\[
M = A^{\top} C A
\]
(this is not a local inverse formula, and it happens only in the special case that the nonzero pattern of $M$ is subordinate to the same chordal graph associated with $M^{-1}$).
The matrix $A$ is
\[
A =
\left(
\begin{tabular}{ccccccc}
      1  & 0     & 0     & 0     & 0     & 0 \\
      0     & 1     & 0     & 0     & 0     & 0  \\
      0     & 0     & 1     & 0     & 0     & 0 \\
     0     & 1     & 0     & 0     & 0     & 0 \\
      0     & 0     & 1     & 0     & 0     & 0 \\
      0     & 0     & 0     & 1     & 0     & 0 \\
      0     & 0     & 1     & 0     & 0     & 0 \\
      0     & 0     & 0     & 1     & 0     & 0 \\
      0     & 0     & 0     & 0     & 1     & 0 \\
      0     & 0     & 0     & 1     & 0     & 0 \\
      0     & 0     & 0     & 0     & 1     & 0 \\
      0     & 0     & 0     & 0     & 0     & 1 \\
      0     & 1     & 0     & 0     & 0     & 0 \\
      0     & 0     & 1     & 0     & 0     & 0 \\
      0     & 0     & 1     & 0     & 0     & 0 \\
      0     & 0     & 0     & 1     & 0     & 0 \\
      0     & 0     & 0     & 1     & 0     & 0 \\
      0     & 0     & 0     & 0     & 1     & 0
\end{tabular}
\right).
\]
The matrix $C$ is block diagonal with blocks, in this order,
\[
C_1=
\left(
\begin{tabular}{ccc}
   -0.8660 &   0.5000  &       0 \\
    0.4830  &  0.8365  &  0.2241\\
    0.1294  &  0.2241  & -0.8365
\end{tabular}
\right),
\;
C_2=
\left(
\begin{tabular}{ccc}
    0.8365 &   0.2241 &  -0.1294 \\
    0.2241 &  -0.8365 &   0.4830 \\
         0  &  0.4830  &  0.8365
\end{tabular}
\right),
\]
\[
C_3=
\left(
\begin{tabular}{ccc}
   -0.8365 &   0.4830 &        0 \\
    0.4830  &  0.8365  &  0.2241 \\
    0.1294  &  0.2241 &  -0.8365
\end{tabular}
\right),
\;
C_4=
\left(
\begin{tabular}{ccc}
    0.8365 &    0.2241&   -0.1294 \\
    0.2241  & -0.8365 &   0.4830 \\
         0   & 0.5000 &   0.8660
\end{tabular}
\right),
\]
\[
-C_5=
\left(
\begin{tabular}{cc}
   -0.8365  & -0.2241\\
   -0.2241 &   0.8365
\end{tabular}
\right), \;
-C_6=
\left(
\begin{tabular}{cc}
    0.8365 &  -0.4830 \\
   -0.4830  & -0.8365
\end{tabular}
\right),
\]
\[
-C_7=
\left(
\begin{tabular}{cc}
   -0.8365 &  -0.2241\\
   -0.2241 &   0.8365
\end{tabular}
\right).
\]
For the special class of matrices for which the local inverse formula applies in both directions,  and analogous to the way a block diagonal $C$ is defined from that part of $M$ corresponding to the chordal graph of $M_0$, we could also define a block diagonal matrix $D$ from that part of $M^{-1}$ corresponding to the same chordal graph.
Then
\begin{equation}
A^{\top} D^{-1} A = M
\label{eq:local:inverse:formula:backward}
\end{equation}
(compared to \eqref{eq:local:inverse:formula:forward}, here \eqref{eq:local:inverse:formula:backward} is the local inverse formula in the \textit{opposite} direction, by assembling $M$ from inverses of blocks and overlaps in $M^{-1}$)
\textit{and}
\[
M^{-1} = A^{\top} D A.
\]
In this example $D$ is the same as $C^\top$, but there are other examples for which the local inverse formula applies in both directions where $D \ne C^\top$.

\paragraph{Acknowledgement}
The authors gratefully acknowledge a grant from The Mathworks that made this work possible.


\begin{thebibliography}{10}
\providecommand{\url}[1]{{#1}}
\providecommand{\urlprefix}{URL }
\expandafter\ifx\csname urlstyle\endcsname\relax
  \providecommand{\doi}[1]{DOI~\discretionary{}{}{}#1}\else
  \providecommand{\doi}{DOI~\discretionary{}{}{}\begingroup
  \urlstyle{rm}\Url}\fi

\bibitem{PeterBartlett2009}
Bartlett, P.: Undirected graphical models: Chordal graphs, decomposable graphs,
  junction trees, and factorizations (2009).
\newblock
  \urlprefix\url{https://people.eecs.berkeley.edu/~bartlett/courses/2009fall-cs281a/}

\bibitem{BlairPeytonReport1993}
Blair, J.R.S., Peyton, B.: Graph Theory and Sparse Matrix Computation,
  \emph{The IMA Volumes in Mathematics and its Applications}, vol.~56, chap. An
  Introduction to Chordal Graphs and Clique Trees, pp. 1--29.
\newblock Springer (1993)

\bibitem{DaubechiesBook}
Daubechies, I.: Ten Lectures on Wavelets.
\newblock Society for Industrial and Applied Mathematics (1992).
\newblock \doi{10.1137/1.9781611970104}

\bibitem{Dempster72}
Dempster, A.P.: Covariance selection.
\newblock Biometrics \textbf{28}, 157--175 (1972)

\bibitem{DymGohberg1981}
Dym, H., Gohberg, I.: Extensions of band matrices with band inverses.
\newblock Linear Algebra and Its Applications \textbf{36}, 1--24 (1981)

\bibitem{EidelmanGohbergHaimovici2013}
Eidelman, Y., Gohberg, I., Haimovici, I.: Separable Type Representations of
  Matrices and Fast Algorithms: Volume 1.
\newblock Springer (2013)

\bibitem{friedman2008sparse}
Friedman, J., Hastie, T., Tibshirani, R.: Sparse inverse covariance estimation
  with the graphical lasso.
\newblock Biostatistics \textbf{9}(3), 432--441 (2008).
\newblock \doi{10.1093/biostatistics/kxm045}

\bibitem{JohnsonProceedings1989}
Johnson, C.R.: {Matrix Completion Problems: A Survey}.
\newblock In: C.R. Johnson (ed.) Matrix Theory and Applications, pp. 69--87.
  American Mathematical Society (1989)

\bibitem{JohnsonLocalInverse1998}
Johnson, C.R., Lundquist, M.: Local inversion of matrices with sparse inverses.
\newblock Linear Algebra and Its Applications \textbf{277}, 33--39 (1998)

\bibitem{KollerFriedman2009}
Koller, D., Friedman, N.: Probabilistic Graphical Models: Principles and
  Techniques.
\newblock MIT Press (2009)

\bibitem{GraphicalModelsLauritzenBook}
Lauritzen, S.: Graphical Models.
\newblock Oxford University Press (1996)

\bibitem{MallatWaveletBook}
Mallat, S.: A Wavelet Tour of Signal Processing.
\newblock Academic Press (1998)

\bibitem{ravikumar2011high}
Ravikumar, P., Wainwright, M.J., Raskutti, G., Yu, B., et~al.: High-dimensional
  covariance estimation by minimizing l1-penalized log-determinant divergence.
\newblock Electronic Journal of Statistics \textbf{5}, 935--980 (2011).
\newblock \doi{doi:10.1214/11-EJS631}

\bibitem{DonRose70}
Rose, D.: Triangulated graphs and the elimination process.
\newblock Journal of Mathematical Analysis and Applications \textbf{32},
  597--609 (1970)

\bibitem{SpeedKiiveri1986}
Speed, T.P., Kiiveri, H.T.: Gaussian markov distributions over finite graphs.
\newblock The Annals of Statistics \textbf{14}(1), 138--150 (1986)

\bibitem{Strang:2010aa}
Strang, G.: Fast transforms: Banded matrices with banded inverses.
\newblock Proc Natl Acad Sci U S A \textbf{107}(28), 12,413--6 (2010).
\newblock \doi{10.1073/pnas.1005493107}

\bibitem{Str16}
Strang, G.: Introduction to {L}inear {A}lgebra.
\newblock Wellesley-Cambridge Press (2016)

\bibitem{StrangWaveletBook}
Strang, G., Nguyen, T.: Wavelets and Filter Banks.
\newblock Wellesley-Cambridge Press (1996)

\bibitem{StrangNguyenInterplay2004}
Strang, G., Nguyen, T.: The interplay of ranks of submatrices.
\newblock SIAM Review \textbf{46}(4), 637--646 (2004).
\newblock \doi{10.1137/S0036144503434381}

\bibitem{SemiSeparableBook}
Vandebril, R., van Barel, M., Mastronardi, N.: Matrix Computations and
  Semiseparable Matrices, vol.~1.
\newblock Johns Hopkins (2007)

\end{thebibliography}
\end{document}